\documentclass[12pt]{article}
\usepackage{latexsym}
\usepackage{times}
\usepackage{graphicx,color}
\usepackage{amsmath,amssymb, amsfonts}
\catcode`@=11
\renewcommand{\@begintheorem}[2]{\it \trivlist            
      \item[\hskip \labelsep{\bf #1\ #2{\rm :}}]}         
\renewcommand{\@opargbegintheorem}[3]{\it \trivlist       
      \item[\hskip \labelsep{\bf #1\ #2\ {\rm (#3)\/:}}]}
\def\@sect#1#2#3#4#5#6[#7]#8{\ifnum #2>\c@secnumdepth
     \def\@svsec{}\else 
     \refstepcounter{#1}\edef\@svsec{\csname the#1\endcsname{.}\hskip 1em }\fi
     \@tempskipa #5\relax
      \ifdim \@tempskipa>\z@ 
        \begingroup #6\relax
          \@hangfrom{\hskip #3\relax\@svsec}{\interlinepenalty \@M #8\par}
        \endgroup
       \csname #1mark\endcsname{#7}\addcontentsline
         {toc}{#1}{\ifnum #2>\c@secnumdepth \else
                      \protect\numberline{\csname the#1\endcsname}\fi
                    #7}\else
        \def\@svsechd{#6\hskip #3\@svsec #8\csname #1mark\endcsname
                      {#7}\addcontentsline
                           {toc}{#1}{\ifnum #2>\c@secnumdepth \else
                             \protect\numberline{\csname the#1\endcsname}\fi
                       #7}}\fi
     \@xsect{#5}}

\@addtoreset{equation}{section}
\newcommand{\Delete}[1]{}
\newcommand{\pend}{\hspace*{\fill} $\Box$}

\newtheorem{lemma}{Lemma}[section]
\newtheorem{theorem}[lemma]{Theorem}

\newtheorem{remark}{Remark}

\begin{document}

\title{A Note on Ordinal Submodularity}

\author{{Satoru Fujishige}\footnote{ 
Research Institute for Mathematical Sciences,
Kyoto University, Kyoto 606-8502, Japan.
\texttt{{\rm Email:} fujishig@kurims.kyoto-u.ac.jp}}\ \ \;
and\ \  \;
{Ryuhei Mizutani}\footnote{
Faculty of Science and Technology, Keio University, Kanagawa  
223-8522, Japan. 
\texttt{{\rm Email:} mizutani@math.keio.ac.jp}}}

\date{February 17, 2026}

\maketitle

\begin{abstract}

Notions of ordinal submodularity/supermodularity have been introduced 
and studied in the literature.
We consider several classes of ordinally submodular functions defined 
on finite Boolean lattices 
and give characterizations of the set of minimizers of ordinally 
submodular functions.

\end{abstract}

\noindent
{\bf Keywords}: Ordinal submodularity, ordinal submodular minimization,
quasisubmodular functions

\section{Introduction}\label{sec:1}
Notions of ordinal submodularity/supermodularity have been introduced 
and studied in the literature (see, e.g., 
\cite{ChambersEchenique2008,ChambersEchenique2009,MilgromShannon1994,MurotaShioura2003,Topkis98}).
For simplicity we consider functions on finite Boolean lattices and investigate 
a characterization of the set of minimizers of ordinally submodular
functions.
The following arguments can basically be extended to the case where functions 
are defined on the integer vector lattice $\mathbb{Z}^E$ with lattice 
operations join $\vee$ and meet $\wedge$ defined by 
$x\vee y=(\max\{x(e),y(e)\}\mid e\in E)$ and 
$x\wedge y=(\min\{x(e),y(e)\}\mid e\in E)$ for all $x,y\in\mathbb{Z}^E$
and even for more general lattices.

In Section~\ref{sec:2} we examine various notions of ordinally submodular 
functions. Specifically we define {\bf (Q1)}-, {\bf (Q2)}-, {\bf (Q3)}-, and 
{\bf (Q4)}-submodularity and discuss the relationship among these notions.
We show characterizations (Lemma~\ref{lem:1a}, Theorem~\ref{th:1})  
of minimizers of {\bf (Q1)}-submodular functions 
in Section~\ref{sec:3.1} and also that of {\bf (Q4)}-submodular 
functions under some additional condition in Section~\ref{sec:3.2}.
Moreover, we consider a very recent result obtained by 
Mizutani~\cite{Mizutani2026} from the point of view of ordinal submodularity,
in Section~\ref{sec:concl}.

\section{Ordinal Notions of Submodularity}
\label{sec:2}

Let $E$ be a nonempty finite set and consider a set function 
$f: 2^E\to(\mathbb{P},\le)$, where $(\mathbb{P},\le)$ is a totally 
ordered set ($\mathbb{P}$ being not necessarily the set of reals or rationals). 
Note that by considering the dual ordered set $(\mathbb{P},\le^*)$ of 
$(\mathbb{P},\le)$ (i.e., 
$\forall a,b\in\mathbb{P}: a\le^* b\Leftrightarrow b\le a$) ordinal notions of 
submodularity associated with $(\mathbb{P},\le)$ are those of supermodularity
associated with $(\mathbb{P},\le^*)$.

For any $X,Y\in 2^E$ we have several ordinal conditions in the literature 
as follows.
\begin{itemize}
  \item[{\bf (Q1)}] 
      $f(X)\le f(X\cap Y) \Longrightarrow f(X\cup Y)\le f(Y)$\;.
  \item[{\bf (Q2)}] $f(X)<f(X\cap Y) \Longrightarrow f(X\cup Y)< f(Y)$\;.
  \item[{\bf (Q3)}]$f(X)<f(X\cap Y) \Longrightarrow f(X\cup Y)\le f(Y)$\;.
  \item[{\bf (Q4)}] $\max\{f(X),f(Y)\}\ge\min\{f(X\cup Y),f(X\cap Y)\}$\;.
\end{itemize}
Various notions of ordinal submodularity are defined and termed as follows.
\begin{enumerate}
\item When Conditions {\bf (Q1)} and {\bf (Q2)} hold for all $X,Y\in 2^E$,
 $f$ is called {\it quasisubmodular} in \cite{MilgromShannon1994}. 
(Such an $f$ is also called {\it semistrictly quasisubmodular} 
in \cite{MurotaShioura2003}.)
\item When Condition {\bf (Q2)} holds for all $X,Y\in 2^E$,
 $f$ is called {\it weakly quasisubmodular} in \cite{ChambersEchenique2008}.
\item  When Condition {\bf (Q3)} holds for all $X,Y\in 2^E$,
 $f$ is called {\it quasisubmodular} in \cite{MurotaShioura2003} 
(but we do not adopt this terminology in the present note.) 
\item Condition {\bf (Q4)} appears as (QSB\_w) in \cite{MurotaShioura2003}.
\end{enumerate}

To avoid any confusion, for each $i=1,2,3,4$ we call $f$ a  
{\bf (Qi)}-{\it submodular function} if it satisfies Condition {\bf (Qi)} 
for all $X,Y\in 2^E$.

It should be noted that {\bf (Q1)}, {\bf (Q2)}, and {\bf (Q3)} are
equivalently described as follows.
\begin{itemize}
  \item[{\bf (Q1)}] 
      $f(X)> f(X\cap Y)$ \ {\rm or}\  $f(X\cup Y)\le f(Y)$\;.
  \item[{\bf (Q2)}] $f(X)\ge f(X\cap Y)$ \ {\rm or}\  $f(X\cup Y)< f(Y)$\;.
  \item[{\bf (Q3)}]$f(X)\ge f(X\cap Y)$ \ {\rm or}\  $f(X\cup Y)\le f(Y)$\;.
\end{itemize}

\begin{remark}{\rm 
Conditions {\bf (Q3)} and {\bf (Q4)} are self-dual for the 
Boolean lattice $2^E$ with respect to the interchanging of the roles
of set union $\cup$ and intersection $\cap$. But this is not the case 
for Condition {\bf (Q1)} or {\bf (Q2)}. 
Actually {\bf (Q2)} is in a dual form of {\bf (Q1)}, so that 
imposing both {\bf (Q1)} and {\bf (Q2)}, we get a self-dual condition
(for quasisubmodularity \cite{MilgromShannon1994}).
Note that {\sl ordinary} submodular functions (\cite{Fuji05}) are defined 
in a self-dual manner.
\pend}
\end{remark}

\begin{remark}{\rm 
Any quasisubmodular function is {\bf (Q1)}-submodular and 
{\bf (Q2)}-submodular. Any {\bf (Q1)}-submodular or 
 {\bf (Q2)}-submodular function is 
{\bf (Q3)}-submodular. 
Functions that satisfy Condition {\bf (Q1)} or {\bf (Q2)} for every
$X,Y \in 2^E$ are exactly {\bf (Q3)}-submodular functions.
\pend}
\end{remark}

\begin{figure}[h]
  \centering
  \includegraphics[width=9cm,height=7cm,clip]{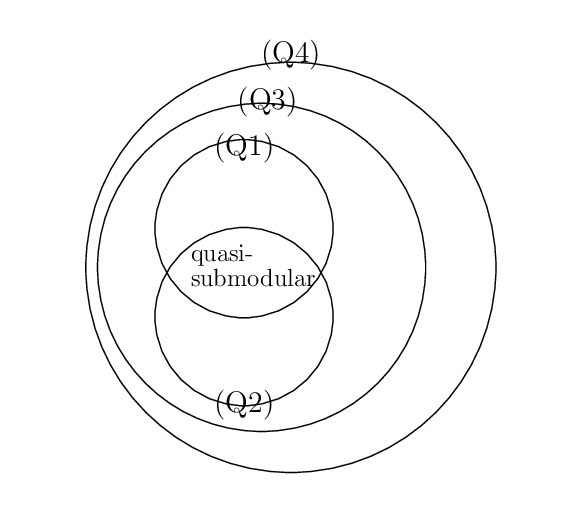}
\caption{Relationship among {\bf (Qi)} for $i=1,2,3,4$.}
\label{fig:1}
\end{figure}


Moreover, we have the following lemma.

\begin{lemma}\label{lem:1}
Any {\bf (Q3)}-submodular function is {\bf (Q4)}-submodular.
\medskip\\
{\rm (Proof)
Let $f: 2^E\to(\mathbb{P},\le)$ be any {\bf (Q3)}-submodular function and 
consider arbitrary $X,Y\in 2^E$.  

If $f(X)\ge f(X\cap Y)$, then 
we have 
\begin{equation} 
\max\{f(X),f(Y)\}\ge f(X)\ge f(X\cap Y)\ge 
\min\{f(X\cup Y),f(X\cap Y)\}\;.
\end{equation}
Also, if $f(X)< f(X\cap Y)$,\; from {\bf (Q3)}
we have $f(X\cup Y)\le f(Y)$.
Hence  we have 
\begin{equation} 
 \max\{f(X),f(Y)\}\ge f(Y)\ge f(X\cup Y)\ge 
\min\{f(X\cup Y),f(X\cap Y)\}\;.
\end{equation}
It follows that the max-min relation of {\bf (Q4)} 
holds for all $X, Y\in 2^E$.
\pend}
\end{lemma}

\begin{remark}\label{rem:3}{\rm 
Suppose for example that some $X,Y\in 2^E$ satisfy 
\begin{equation}
 f(X)<f(X\cap Y)< f(Y)<f(X\cup Y)\;.
\end{equation}
Then the max-min relation of {\bf (Q4)} holds for the present $X, Y$
but {\bf (Q3)} does not hold.
Hence {\bf (Q4)}-submodularity is strictly weaker than 
{\bf (Q3)}-submodularity. See Figure~\ref{fig:1}.
\pend}
\end{remark}

\section{Minimizers of Ordinally Submodular Functions}\label{sec:3}

For any $X,Y\in 2^E$ with $X\subseteq Y$ define 
$[X,Y]=\{Z\in 2^E\mid X\subseteq Z\subseteq Y\}$.

\subsection{{\bf (Q1)}-submodular functions}\label{sec:3.1}

In this subsection we consider {\bf (Q1)}-submodularity, which 
is weaker than  quasisubmodularity satisfying both {\bf (Q1)} and 
{\bf (Q2)}.

We have the following lemma.

\begin{lemma}\label{lem:1a}
Let $f: 2^E\to(\mathbb{P},\le)$ be a {\bf (Q1)}-submodular function. 
For any $X\in 2^E$, 
if $X$ is a minimizer of $f$ over $[\emptyset,X]$,
then there exists a global minimizer $Z^*$ of $f$ over $2^E$ such that
$Z^*\in[X,E]$.
\medskip\\
{\rm (Proof)
Let $Z^*$ be a global minimizer of $f$ over $2^E$. Given a set $X\in 2^E$, 
suppose that $X$ is a minimizer of $f$ over $[\emptyset,X]$.  
Since $Z^*\cap X\in[\emptyset,X]$, it follows from the assumption on $X$
that we have $f(X)\le f(Z^*\cap X)$.
Then because of {\bf (Q1)} we have $f(Z^*\cup X)\le f(Z^*)$. 
This means that $Z^*\cup X$ belonging to $[X,E]$ is a global minimizer of 
$f$ over $2^E$.
\pend}
\end{lemma}

The proof of the present lemma is a direct adaptation of that of 
\cite[Lemma~1]{FujiIwata2002} for ordinary submodular functions.

{F}rom Lemma~\ref{lem:1a} we get the following theorem.

\begin{theorem}\label{th:1}
Let $f: 2^E\to(\mathbb{P},\le)$ be a {\bf (Q1)}-submodular function. 
If $X^*\in 2^E$ is a minimizer of $f$ over $[\emptyset,X^*]\cup[X^*,E]$,
then $X^*$ is also a minimizer of $f$ over $2^E$.
\medskip\\
{\rm (Proof) 
Suppose that 
$X^*\in 2^E$ is a minimizer of $f$ over $[\emptyset,X^*]\cup[X^*,E]$.
Then it follows from Lemma~\ref{lem:1a} that there exists a global 
minimizer $Z^*$ of $f$ over $2^E$ that belongs to $[X^*,E]$.
Hence we must have $f(X^*)=f(Z^*)$.
\pend}
\end{theorem}

\begin{remark}{\rm
Theorem~\ref{th:1} gives a local characterization of global minimizers 
of a {\bf (Q1)}-submodular function $f$. 
Similarly as Lemma~\ref{lem:1a} and Theorem~\ref{th:1} 
a related local characterization of global minimizers of {\sl ordinary}  
submodular functions is known (see \cite[Theorem~7.2]{Fuji05}). 
Recall that {\sl ordinary} submodular functions are quasisubmodular functions
taking values on a totally ordered additive group $\mathbb{P}$.
Also note that Theorem~\ref{th:1} for a quasisubmodular function is also 
shown in \cite[Corollary~6.3]{MurotaShioura2003}, where quasisubmodularity 
is stronger than {\bf (Q1)}-submodularity. 
\pend}
\end{remark}

\begin{remark}{\rm
If $f$ is quasisubmodular, i.e., $f$ satisfies both {\bf (Q1)} and 
{\bf (Q2)} for all $X,Y\in 2^E$, then putting 
$\mathcal{D}_f={\rm Arg}\min\{f(X)\mid X\in 2^E\}$, for all 
$X,Y\in \mathcal{D}_f$ we have 
$X\cup Y, X \cap Y\in \mathcal{D}_f$. 
For, if $X\cup Y\not\in\mathcal{D}_f$, then 
$f(X\cup Y)>f(Y)(=f(X))$, which implies from {\bf (Q1)} $f(X\cap Y)<f(X)$, 
a contradiction. Similarly we get a contradiction if $f(X\cap Y)>f(X)(=f(Y))$, 
due to {\bf (Q2)}. (Also see \cite[Theorem~5.8]{MurotaShioura2003}.)
\pend}
\end{remark}

\begin{remark}{\rm
We can show a dual form of Lemma~\ref{lem:1a} for {\bf (Q2)}-submodular
functions by interchanging lattice operations $\cup$ and $\cap$, which 
also leads us to Theorem~\ref{th:1} for {\bf (Q2)}-submodular
functions. 
\pend}
\end{remark}

\subsection{{\bf (Q4)}-submodular functions}\label{sec:3.2}

Let us consider a {\bf (Q4)}-submodular function $f: 2^E\to(\mathbb{P},\le)$.  
We assume that 
\begin{itemize}
\item[{\bf $(*)$}] $f$ defines a linear order on the 
power set $2^E$, i.e., $f(X)\neq f(Y)$ for all distinct $X, Y\in 2^E$.
\end{itemize}
Under Assumption $(*)$ we have the following theorem for 
{\bf (Q4)}-submodular functions similar to 
Theorem~\ref{th:1} for {\bf (Q1)}-submodular functions.

\begin{theorem}\label{th:2}
Let $f: 2^E\to(\mathbb{P},\le)$ be a {\bf (Q4)}-submodular function 
 satisfying $(*)$. 
If $X^*\in 2^E$ is a minimizer of $f$ over $[\emptyset,X^*]\cup[X^*,E]$,
then $X^*$ is also a minimizer of $f$ over $2^E$.
\medskip\\
{\rm (Proof)
Let $Z^*$ be a minimizer of $f$ over $2^E$ (a unique global minimizer of $f$
due to Assumption $(*)$).
 Then it follows from {\bf (Q4)} and the definitions of $X^*$ and $Z^*$ 
that we have 
\begin{equation}\label{eq:2a2}
 f(X^*)=\max\{f(X^*),f(Z^*)\}\ge\min\{f(X^*\cap Z^*),f(X^*\cup Z^*)\}
 \ge f(X^*)\;.
\end{equation}
Hence we have 
\begin{equation}\label{eq:2b2}
  f(X^*)=\min\{f(X^*\cup Z^*),f(X^*\cap Z^*)\}\;.
\end{equation}
It follows from (\ref{eq:2b2}) that
we have $f(X^*)=f(X^*\cup Z^*)$ or $f(X^*)=f(X^*\cap Z^*)$.
Because of Assumption $(*)$ we have $X^*=X^*\cup Z^*$ or 
$X^*=X^*\cap Z^*$, which implies $Z^*\in [\emptyset,X^*]\cup[X^*,E]$.
Because of Assumption $(*)$ and the definitions of $X^*$ and $Z^*$
we must have $Z^*=X^*$.
\pend}
\end{theorem}

It should be noted that a {\bf (Q4)}-submodular function with property $(*)$
is not {\bf (Q3)}-submodular in general (see Remark~\ref{rem:3}).

\section{Discussions}\label{sec:concl}

Let us consider another ordinal function $f: 2^E\to(\mathbb{P},\le)$ 
that satisfies for all $X,Y\in 2^E$ 
\begin{itemize}
\item[{\bf (Qh)}] if $f(X)=f(Y)$, then $f(X\cup Y)=f(X\cap Y)=f(X)(=f(Y))$ or
$f(X\cup Y)<f(X)(=f(Y))$ or $f(X\cap Y)<f(X)(=f(Y))$.
\end{itemize}
We call such a function a {\it {\bf (Qh)}-submodular function}.
The {\bf (Qh)}-submodularity is closely related to the concept of 
{\sl hierarchical lattice} which has very recently been introduced by
Mizutani \cite{Mizutani2026}.
Note that every quasisubmodular function is a {\bf (Qh)}-submodular 
function but the converse is not true in general. 

For a {\bf (Qh)}-submodular function $f: 2^E\to(\mathbb{P},\le)$ 
let the distinct values of $f(X)$ for all $X\in 2^E$ be given by
$\mu_1<\mu_2<\cdots<\mu_p$. For each $i=1,2,\cdots,p$ define 
$\mathcal{F}_i=\{X\in 2^E\mid f(X)\le \mu_i\}$. Then we have a sequence
of families:
\begin{equation}\label{eq:4-1}
\emptyset=\mathcal{F}_0\subset\mathcal{F}_1\subset\cdots
\subset\mathcal{F}_p=2^E.
\end{equation}
For any $k\in\{1,\cdots,p\}$ $\mathcal{F}_k$ is exactly the family called 
a {\it $k$-hierarchical lattice} in \cite{Mizutani2026}. 
(Note that when we are given a $k$-hierarchical lattice in the original 
sense of the definition in \cite{Mizutani2026}, we can construct 
a {\bf (Qh)}-submodular function $f$ such that its induced family 
$\mathcal{F}_k$ coincides with the given $k$-hierarchical lattice.) 
Moreover, Mizutani \cite{Mizutani2026} considered the problem of minimizing 
an {\sl ordinary} submodular function $\varphi: 2^E\to\mathbb{R}$ 
(the set of reals) subject to 
$X\in 2^E\setminus \mathcal{F}_k$. That is, 
for a fixed $k\in\{1,\cdots,p$$-$$1\}$
\begin{equation}\label{eq:4-2}
 {\rm Minimize}\ \ \varphi(X)\ \ {\rm subject\ to }\ \ 
    f(X)> \mu_k\;.
\end{equation}
The arguments made in \cite{Mizutani2026} except for those related to 
the algorithm complexity can basically be extended to the case where
the functions $\varphi$ and $f$ are quasisubmodular (or each being a 
slightly more general ordinal function).

\section*{Acknowledgements}

This work was supported by JST ERATO Grant Number JPMJER2301, Japan.
S.~Fuji\-shige's research was supported by JSPS KAKENHI Grant Numbers 
JP22K11922 and by the Research Institute for Mathematical Sciences, 
an International Joint Usage/Research Center located in Kyoto University.


\begin{thebibliography}{99}

\bibitem{ChambersEchenique2008} C. P. Chambers and F. Echenique: 
Ordinal notions of submodularity. {\it Journal of Mathematical Economics}
{\bf 44} (2008) 1243--1245.

\bibitem{ChambersEchenique2009} C. P. Chambers and F. Echenique: 
Supermodularity and preferences. {\it Journal of Economic Theory}
{\bf 144} (2009) 1004--1014.

\bibitem{Fuji05} S. Fujishige: {\it Submodular Functions and 
Optimization}, North-Holland, 1991 (Second Edition, Elsevier, 2005).

\bibitem{FujiIwata2002} S. Fujishige and S. Iwata: A descent method for
submodular function minimization. {\it Mathematical Programming}, Ser.~A 
{\bf 92} (2002) 387--390.

\bibitem{MilgromShannon1994} P. Milgrom and C. Shannon: Monotone comparative 
statics. {\it Econometrica} {\bf 62} (1994) 157--180.

\bibitem{Mizutani2026} R. Mizutani: Minimizing submodular functions over 
hierarchical families. arXiv:2601.14805v1 [math.CO] 21 Jan 2026.

\bibitem{MurotaShioura2003} K. Murota and A. Shioura: Quasi M-convex and 
L-convex functions---quasiconvexity in discrete optimization.
{\it Discrete Applied Mathematics} {\bf 131} (2003) 467--494.

\bibitem{Topkis98} D. M. Topkis: {\it Supermodularity and Complementarity} 
(Princeton University Press, 1998).


\end{thebibliography}
\end{document}